\documentclass[a4paper,10pt]{scrartcl}
\pdfoutput=1
\usepackage{ifpdf}
\usepackage[english]{babel}
\usepackage[utf8]{inputenc}
\usepackage[T1]{fontenc}
\usepackage{lmodern}
\usepackage{enumerate}
\usepackage{geometry}
\usepackage{amsfonts,amsmath,amssymb,amsopn}
\usepackage{amsthm}
\usepackage{mathtools}
\usepackage{stmaryrd}
\usepackage{nicefrac}
\usepackage{exscale}
\usepackage{empheq}
\usepackage{ mathrsfs }

\usepackage[numbers,square]{natbib}
\usepackage{url} 
\usepackage[colorlinks=true,pdfpagelabels,unicode]{hyperref}
\hypersetup{linkcolor=blue, urlcolor=blue, citecolor=blue}
\usepackage[dvipsnames,svgnames,table]{xcolor}
\usepackage{orcidlink} 

\allowdisplaybreaks 

\theoremstyle{plain}

\newtheoremstyle{theo}
	{3pt} 
	{3pt} 
	{\itshape} 
	{} 
		{\bfseries} 
	{\\} 
	{ } 
	{\thmname{#1}\thmnumber{ #2.}\thmnote{ - #3}} 
\theoremstyle{theo}

\newtheorem{definition}{Definition}[section]
\newtheorem{lemma}[definition]{Lemma}
\newtheorem{theorem}[definition]{Theorem}

\newtheorem{proposition}[definition]{Proposition}

\newenvironment{bew}{\begin{proof}[\bfseries Proof:]}{\end{proof}}

\newtheoremstyle{remark}
	{3pt} 
	{3pt} 
	{} 
	{} 
		{\bfseries} 
	{} 
	{ } 
	{\thmname{#1}\thmnumber{ #2.}\thmnote{ - #3}} 
\theoremstyle{remark}
\newtheorem{remark}[definition]{Remark}

\DeclareMathOperator{\bomega}{\overline{\Omega}}
\DeclareMathOperator{\romega}{\partial\Omega}
\DeclareMathOperator{\intd}{d\!}

\newcommand{\Tm}{T_{max}}

\newcommand{\Lo}[1][1]{L^{#1}(\Omega)} 
\newcommand{\W}[1][1,2]{W^{#1}(\Omega)}

\newcommand{\LSp}[2]{L^{#1\;\!}\!\left(#2\right)} 

\newcommand{\CSp}[2]{C^{#1}\!\left(#2\right)}
\newcommand{\CSpb}[2]{C^{#1}\!\big(#2\big)}
\newcommand{\CSpnl}[2]{C^{#1}\!\,(#2)} 

\newcommand{\R}{\mathbb{R}}

\newcommand{\dimN}{N}

\newcommand{\matr}[1]{\begin{pmatrix}#1\end{pmatrix}}
\newcommand{\smatr}[1]{\left(\begin{smallmatrix}#1\end{smallmatrix}\right)}

\makeatletter
\def\@fnsymbol#1{\ensuremath{\ifcase#1\or *\or \ddagger\or \#\or
   \mathsection\or \mathparagraph\or \|\or **\or \dagger\dagger
   \or \ddagger\ddagger \else\@ctrerr\fi}}
\def\@fnsymbol#1{\ensuremath{\ifcase#1\or *\or \ddagger\or \#\or
   \mathsection\or \mathparagraph\or \|\or **\or \dagger\dagger
   \or \ddagger\ddagger \else\@ctrerr\fi}}
   
 \def\@lbibitem[#1]#2#3{%
  \if\relax\@extra@b@citeb\relax\else
    \@ifundefined{br@#2\@extra@b@citeb}{}{%
     \@namedef{br@#2}{\@nameuse{br@#2\@extra@b@citeb}}%
    }%
  \fi
  \@ifundefined{b@#2\@extra@b@citeb}{%
   \def\NAT@num{}%
  }{%
   \NAT@parse{#2}%
  }%
  \def\NAT@tmp{#1}%
  \expandafter\let\expandafter\bibitemOpen\csname NAT@b@open@#2\endcsname
  \expandafter\let\expandafter\bibitemShut\csname NAT@b@shut@#2\endcsname
  \@ifnum{\NAT@merge>\@ne}{%
   \NAT@bibitem@first@sw{%
    \@firstoftwo
   }{%
    \@ifundefined{NAT@b*@#2}{%
     \@firstoftwo
    }{%
     \expandafter\def\expandafter\NAT@num\expandafter{\the\c@NAT@ctr}%
     \@secondoftwo
    }%
   }%
  }{%
   \@firstoftwo
  }%
  {%
   \global\advance\c@NAT@ctr\@ne
   \@ifx{\NAT@tmp\@empty}{\@firstoftwo}{%
    \@secondoftwo
   }%
   {%
    \expandafter\def\expandafter\NAT@num\expandafter{\the\c@NAT@ctr}%
    \global\NAT@stdbsttrue
   }{}%
   \bibitem@fin
   \item[\href{#3}{\hfil\NAT@anchor{#2}{\NAT@num}}]
   \global\let\NAT@bibitem@first@sw\@secondoftwo
   \NAT@bibitem@init
  }%
  {%
   \NAT@anchor{#2}{}%
   \NAT@bibitem@cont
   \bibitem@fin
  }%
  \@ifx{\NAT@tmp\@empty}{%
    \NAT@wrout{\the\c@NAT@ctr}{}{}{}{#2}%
  }{%
    \expandafter\NAT@ifcmd\NAT@tmp(@)(@)\@nil{#2}%
  }%
}
\makeatother 

\author{
Tobias Black\footnote{tblack@math.upb.de}\ \orcidlink{0000-0001-9963-0800}\\
{\small Institute of Mathematics, Paderborn University,}\\[-5pt]
{\small 33098 Paderborn, Germany}
}
\title{Absence of dead-core formations\\in chemotaxis systems with degenerate diffusion}
\date{}

\begin{document}
\maketitle
\begin{abstract}
\noindent
{\textbf{Abstract.} In this paper we consider a chemotaxis system with signal consumption and degenerate diffusion of the form
\begin{align*}
\left\lbrace
\begin{array}{r@{}l@{\quad}l}
&u_t=\nabla\cdot\big(D(u)\nabla u-uS(u)\nabla v\big)+f(u,v),\\
&v_t=\Delta v- uv,\\
\end{array}\right.
\end{align*}
in a bounded domain $\Omega\subset\mathbb{R}^{N}$ with smooth boundary subjected to no-flux and homogeneous Neumann boundary conditions. Herein, the diffusion coefficient $D\in C^0([0,\infty))\cap C^2((0,\infty))$ is assumed to satisfy $D(0)=0$, $D(s)>0$ on $(0,\infty)$, $D'(s)\geq 0$ on $(0,\infty)$ and that there are $s_0>0$, $p>1$ and $C_D>0$ such that $$s D'(s)\leq C_D D(s)\quad\text{and}\quad C_D s^{p-1}\leq D(s)\quad\text{for }s\in[0,s_0].$$ The sensitivity function $S\in C^2([0,\infty))$ and the source term $f\in C^{1}([0,\infty)\times[0,\infty))$ are supposed to be nonnegative.\smallskip

We show that for all suitably regular initial data $(u_0,v_0)$ satisfying $u_0\geq \delta_0>0$ and $v_0\not\equiv 0$ there is a time-local classical solution and -- despite the degeneracy at $0$ -- the solution satisfies an extensibility criterion of the form
$$\text{either}\quad T_{max}=\infty,\quad\text{or}\quad\limsup_{t\nearrow T_{max}}\|u(\cdot,t)\|_{L^\infty(\Omega)}=\infty.$$
Moreover, as a by-product of our analysis, we prove that a classical solution on $\Omega\times(0,T)$ obeying $\|u(\cdot,t)\|_{L^\infty(\Omega)}\leq M_u$ for all $t\in(0,T)$ and emanating from initial data $(u_0,v_0)$ as specified above remains strictly positive throughout $\Omega\times(0,T)$, i.e. one can find $\delta_u=\delta_u(T,\delta_0, M_u,\|v_0\|_{W^{1,\infty}(\Omega)})>0$ such that $$u(x,t)\geq\delta_u\quad\text{for all }(x,t)\in\Omega\times(0,T).$$
Together, the results indicate that the formation of a dead-core in these chemotaxis systems with a degenerate diffusion are impossible before the blow-up time.
}\medskip



\end{abstract}

\newpage
\section{Introduction}\label{sec1:intro}
In reaction-diffusion systems dedicated to the theory of chemical engineering it is quite natural to look for areas of the spatial domain, where the reactants have been consumed completely and no further reaction can take place, this sub-region of the domain is called a dead-core (\cite{temkin1975diffusion,skrzypaczDeadcoreNondeadcoreSolutions2020}).
The systems paramount to our investigation, however, have their origin in the more biological framework of chemotaxis equations. Many bacteria populations act in response to chemical concentration gradients in their surroundings, adjusting their movement-scheme in favor or disfavor of higher signal concentration, all the while consuming or producing some of the substance. Despite the disparate modeling backgrounds between these system families, the search for dead-cores may nevertheless provide valuable qualitative insight. A substantial amount of chemotactically active populations have been witnessed to possess the ability to form patterns of quite magnificent complexity by means of cell aggregation. Examples can be found in \emph{Bacillus subtilis} (\cite{tuval2005bacterial,chertockSinkingMergingStationary2012}), \emph{Escherichia coli} (\cite{budrene91}) and \emph{Salmonella typhirurium} (\cite{woodwardSpatiotemporalPatternsGenerated1995}) to mention a few. (See also \cite[Chapter 5]{murray2003} for general notes on bacterial patterns and chemotaxis.)

Visible patterns in the population, however, are only partially described by areas where the self-organized accumulation of cells happens. Regions featuring an absence of cells are most certainly as crucial in the description of a detectable pattern. (See e.g. depletion zone in the dark-field imaging of experiments concerning self-concentration of \emph{B. subtilis} undertaken in \cite{tuval2005bacterial}.) Accordingly, information on the presence or absence of dead-cores for the bacteria provides additional insight on the process of pattern formation.

We will consider a chemotaxis system of the form
\begin{align}\label{eq:CT-equation}
\left\lbrace
\begin{array}{r@{}l@{\quad}l@{\quad}l@{\,}c}
&u_t=\nabla\cdot\big(D(u)\nabla u -uS(u)\nabla v\big)+f(u,v),\ &x\in\Omega,& t\in(0,T),\\
&v_t=\Delta v -uv,\ &x\in\Omega,& t\in(0,T),\\
&\big(D(u)\nabla u-u S(u)\nabla v\big)\cdot\nu=\nabla v\cdot\nu=0, &x\in\romega,& t\in(0,T),\\
&u(\cdot,0)=u_0,\quad v(\cdot,0)=v_0, &x\in\Omega,&
\end{array}\right.
\end{align}
in a smoothly bounded domain $\Omega\subseteq\R^\dimN$. Herein, 
%
%
we assume $D\in\CSp{0}{[0,\infty)}\cap\CSp{2}{(0,\infty)}$ to satisfy
\begin{align}\label{eq:cond-D}
D(0)=0,\qquad D>0\ \text{ on }\ (0,\infty),\qquad D'\geq0\ \text{ on }\ (0,\infty)
\end{align}
and suppose that there are $s_0>0$, $p>1$ and $C_D>0$ such that
\begin{align}\label{eq:cond-D-hoelder}
sD'(s)\leq C_D D(s)\ \text{ for }\ s\in[0,s_0]\qquad\text{and}\qquad C_D s^{p-1}\leq D(s)\text{ on }[0,s_0].
\end{align}

Moreover, we assume
\begin{align*}
&S\in\CSp{2}{[0,\infty)}\ \text{ and }\ f\in\CSp{1}{[0,\infty)\times[0,\infty)}\ \text{are nonnegative},
\end{align*}
and for $M>0$ and $L>0$ we introduce $C_S=C_S(M)>0$ and $C_f=C_f(M,L)>0$  such that
\begin{align}\label{eq:cond-S_f}
S\leq C_S\ \text{ on }\ [0,M]\ \text{ and }\ f\leq C_f\ \text{ on }\ [0,M]\times[0,L].
\end{align}
The initial data are assumed to fulfill
\begin{align}\label{eq:IR}
\begin{cases}
u_0\in\W[1,\infty]\ \text{is positive with }u_0\geq\delta_0>0, \\
v_0\in\W[1,\infty]\ \text{is nonnegative with }v_0\not\equiv 0.
\end{cases}
\end{align}
A solution $(u,v)$ of \eqref{eq:CT-equation} on $\Omega\times(0,T)$ will be called a dead-core solution if its first component vanishes together with its spatial derivatives inside an interior region of $\Omega$ at some time $t\in(0,T)$. The set $\Omega_{0}^t:=\{x\in\Omega\,\vert\,u(x,t)=0\}\subseteq\Omega$ will be called the dead-core of \eqref{eq:CT-equation} at time $t$ (\cite{bandleDiffusionReactionMonotone1984, bandleFormationDeadCore1984,friedmanFreeBoundarySemilinear1984}).

We will first establish that a solution emanating from a strictly positive initial density distribution of the bacterial population will indeed preserve positivity as long as it remains bounded.

\begin{proposition}\label{prop:1}
Assume that $(u_0,v_0)$ satisfy \eqref{eq:IR}. Let $K_{v_0}:=\|v_0\|_{\W[1,\infty]}>0$, $M_u>0$ and $T\in(0,\infty)$. Then, there is $\delta_u=\delta_u(T,\delta_0,M_u,K_{v_0})>0$ such that if $(u,v)$ comprise a classical solution of \eqref{eq:CT-equation} in $\Omega\times(0,T)$ and if
\begin{align}\label{eq:u-linfty-bound}
\|u(\cdot,t)\|_{\Lo[\infty]}\leq M_u \quad\text{for all }t\in(0,T),
\end{align} then
\begin{align*}
u(x,t)\geq \delta_u\quad\text{for all }(x,t)\in\Omega\times\big(0,T\big).
\end{align*}
\end{proposition}

After that, we will see that this preservation of positivity throughout $\Omega\times(0,T)$ also largely influences the extensibility criterion for local solutions. In fact, from the innate boundedness of $\|v(\cdot,t)\|_{\Lo[\infty]}$ -- enforced by the consumption of signal in the second equation -- and Proposition~\ref{prop:1}, we can show that $\limsup_{t\nearrow\Tm}\|u(\cdot,t)\|_{\Lo[\infty]}$ is the sole deciding factor between in whether the classical solution exists globally in time or not.

\begin{theorem}\label{theo:ex-alt}
Assume $(u_0,v_0)$ satisfy \eqref{eq:IR}. Then there exist $\Tm\in(0,\infty]$ and a uniquely determined pair $(u,v)$ of functions
\begin{align*}
u,v\in \bigcap_{q\geq1}\CSp{0}{\bomega\times[0,\Tm);\W[1,q]}\cap\CSp{2,1}{\bomega\times(0,\Tm)}
\end{align*}
solving \eqref{eq:CT-equation} classically in $\Omega\times(0,\Tm)$. Moreover, $u>0$ in $\bomega\times(0,\Tm)$ and
\begin{align}\label{eq:ref-ex-alt}
\text{either}\quad\Tm=\infty,\quad\text{or}\quad\limsup_{t\nearrow\Tm} \big\|u(\cdot,t)\big\|_{\Lo[\infty]}=\infty.
\end{align}
\end{theorem}
Interpreting both of these results with respect to the initially mentioned dead-core formations, we can conclude that in the system \eqref{eq:CT-equation} a dead-core formation cannot occur before the blow-up time, i.e. $\Omega_0^t=\{\}$ for all $t<\Tm$.

\begin{remark}
Presupposing in Proposition~\ref{prop:1} additionally that $\|v(\cdot,t)\|_{\Lo[\infty]}\leq L_v$ for all $t\in(0,T)$, the strict positivity result for $u$ easily extends to a more general second equation, i.e. $v_t=\Delta v+g(u,v)$ with $g\in\CSp{1}{[0,\infty)\times[0,\infty)}$. In this case the dichotomy in Theorem~\ref{theo:ex-alt} would take the form: $$\text{either}\quad\Tm=\infty,\quad\text{or}\quad\limsup_{t\nearrow\Tm}\big(\big\|u(\cdot,t)\big\|_{\Lo[\infty]}+\big\|v(\cdot,t)\big\|_{\Lo[\infty]}\big)=\infty.$$
For the standard Keller--Segel structure obtained by setting $g(u,v)=-v+u$, one can easily eliminate $\|v(\cdot,t)\|_{\Lo[\infty]}$ from the extensibility criterion by means of a standard semigroup argument and recover \eqref{eq:ref-ex-alt} precisely as in Theorem~\ref{theo:ex-alt}. (See also \cite[Proposition 2.3]{Win24_degen_survey}.)
\end{remark}

\begin{remark}
The second condition on $D$ in \eqref{eq:cond-D} can be replaced by requiring instead the less restricting condition that $\int_0^{s_0}\!\frac{1}{|D(\sigma)|^q}\intd\sigma\leq C_D$ for some $s_0>0$, $q>0$ and $C_D>0$ (\cite[Remark 1.4]{TB23_hoeldertaxis}).
\end{remark}

The main steps of our approach are the following: First, we make use of Amann's results on the solvability of parabolic systems to obtain time-local solutions to \eqref{eq:CT-equation} and a corresponding extensibility criterion. Afterwards, we will make use of semigroup estimates for the Neumann heat semigroup to further eliminate the dependency on $v$ from the existence alternative. Then, assuming $u$ to be bounded, we derive a Hölder bound for $u$, which enables us to establish enhanced regularity properties for $v$ by means of parabolic Schauder-theory. With this, we can then find suitable subsolutions to use in a comparison argument for deriving Proposition~\ref{prop:1}. From which Theorem~\ref{theo:ex-alt} follows by contradiction.

\setcounter{equation}{0} 
\section{Existence of a maximally extended solution}\label{sec2}
To begin, let us affirm the existence of a time-local classical solution and a corresponding extensibility criterion by drawing on well-established general theory (\cite{amann93}). The detailed reasoning herein is along the lines of \cite[Lemma 2.2]{WinLiCPAA22}.
\begin{lemma}\label{lem:exist}
Assume $(u_0,v_0)$ satisfy \eqref{eq:IR}. Then there exist $\Tm\in(0,\infty]$ and a uniquely determined pair $(u,v)$ of functions
\begin{align*}
u,v\in \bigcap_{q\geq1}\CSp{0}{\bomega\times[0,\Tm);\W[1,q]}\cap\CSp{2,1}{\bomega\times(0,\Tm)}
\end{align*}
which constitute a classical solution of \eqref{eq:CT-equation} in $\Omega\times(0,\Tm)$. Moreover, $u>0$ in $\bomega\times(0,\Tm)$ and
\begin{align*}
\text{if }\ \Tm<\infty,\quad\text{then}\quad \limsup_{t\nearrow\Tm}\Big(\|u(\cdot,t)\|_{\Lo[\infty]}+\big\|\tfrac{1}{u(\cdot,t)}\big\|_{\Lo[\infty]}+\|v(\cdot,t)\|_{\Lo[\infty]}\Big)=\infty.
\end{align*}
\end{lemma}

\begin{bew}
Set $D_0:=\R\times(0,\infty)$,
\begin{align*}
A\smatr{\eta\\\xi}:=\matr{1&0\\-\xi S(\xi)&D(\xi)}\quad\text{and}\quad F\smatr{\eta\\\xi}:=\matr{-\eta\xi\\f(\xi,\eta)}.
\end{align*}
For $U\in D_0$ define the operators $$\mathcal{A}\smatr{\eta\\\xi}U:=-\sum_{j=1}^\dimN\partial_j\big(A\smatr{\eta\\\xi}\partial_j U\big)\quad\text{and}\quad \mathcal{B}\smatr{\eta\\\xi}U:=\sum_{j=1}^\dimN\nu_j A\smatr{\eta\\\xi}\partial_j U,$$
where $\nu=(\nu_1,\dots,\nu_\dimN)^T$ denotes the unit outward normal vector at $x\in\romega$. Now, consider the quasilinear problem
\begin{align}\label{eq:gen-eq}
\left\{\begin{array}{r@{\ }l@{\qquad}l}
U_t+\mathcal{A}(U) U&=F(U),&x\in\Omega\times(0,\infty)\\
\mathcal{B}(U)U&=0,&x\in\romega\times(0,\infty),\\
U(\cdot,0)&=U_0,&x\in\Omega.
\end{array}\right.
\end{align}
with $U_0=\smatr{v_0\\u_0}\in\big(\W[1,\infty]\big)^2$. Evidently, $(\mathcal{A},\mathcal{B})$ is lower triangular due to the lower triangular form of $A$ and for each $(\eta,\xi)\in D_0$ the eigenvalues of $A\smatr{\eta\\\xi}$ are positive due to $\xi>0$. Accordingly, $(\mathcal{A},\mathcal{B})$ is normally elliptic and we may employ \cite[Theorem 14.4, Theorem 14.6 and Theorem 15.5]{amann93} to conclude that there are $\Tm\in(0,\infty]$ and a unique $U\in\bigcap_{q\geq1}\CSpb{0}{[0,\Tm);\big(\W[1,q]\big)^2}\cap\big(\CSpb{2,1}{\bomega\times(0,\Tm)}\big)^2$, which solves \eqref{eq:gen-eq} in $\Omega\times(0,\Tm)$ in the classical sense. Moreover, $(v,u)^T:=U$ satisfies $u>0$ in $\bomega\times(0,\Tm)$ and fulfills the extensibility criterion
\[
\text{if}\quad\Tm<\infty,\quad\text{then}\quad\limsup_{t\nearrow\Tm}\Big(\big\|u(\cdot,t)\big\|_{\Lo[\infty]}+\big\|\tfrac{1}{u(\cdot,t)}\big\|_{\Lo[\infty]}+\big\|v(\cdot,t)\big\|_{\Lo[\infty]}\Big)=\infty.
\qedhere\]
\end{bew}

Exploiting the fact that the consumption of the signal chemical in the second equation comes with a innate decay of $\|v(\cdot,t)\|_{\Lo[\infty]}$, we can easily refine the extensibility criterion to be independent of the second component $v$. 

\begin{lemma}\label{lem:v-linfty_bound}
Assume $(u_0,v_0)$ satisfy \eqref{eq:IR} and denote by $(u,v)$ the solution of \eqref{eq:CT-equation} with maximal existence time $\Tm\in(0,\infty]$ provided by Lemma~\ref{lem:exist}. Then,
\begin{align*}
\|v(\cdot,t)\|_{\Lo[\infty]}\leq \|v_0\|_{\Lo[\infty]} \quad\text{for all }t\in[0,\Tm).
\end{align*}
Moreover,
\begin{align}\label{eq:ex-alt-2}
\text{if}\quad\Tm<\infty,\quad\text{then}\quad\limsup_{t\nearrow\Tm}\Big(\big\|u(\cdot,t)\big\|_{\Lo[\infty]}+\big\|\tfrac{1}{u(\cdot,t)}\big\|_{\Lo[\infty]}\Big)=\infty.
\end{align}
\end{lemma}

\begin{bew}
Letting $\bar{v}:=\|v_0\|_{\Lo[\infty]}$ we find from the nonnegativity of $u$ that 
\begin{align*}
\bar{v}_t-\Delta \bar{v}+u\bar{v}=u\bar{v}\geq 0,
\end{align*}
and an application of the parabolic comparison principle entails the desired estimate. The newly specified extensibility criterion is then an evident consequence of this bound.
\end{bew}

The final cultivation of a criterion of the form presented in Theorem~\ref{theo:ex-alt}, however, requires more intricate boundedness properties and their arrangement is the topic of the remaining sections.

\setcounter{equation}{0} 
\section{Hölder regularity of bounded solutions near the maximal existence time}\label{sec3}
In order to prepare the comparison argument undertaken in Section~\ref{sec4}, we aim to establish a bound on $\|\Delta v\|_{\Lo[\infty]}$ for times near the maximal existence time under the assumption that $u$ satisfies \eqref{eq:u-linfty-bound}. To this end, we will perform a step-by-step improvement on the currently obtained regularity properties. 

For the remainder, we fix $(u_0,v_0)$ satisfying \eqref{eq:IR} and denote by $(u,v)$ the classical solution of \eqref{eq:CT-equation} in $\Omega\times(0,T)$. We moreover assume $T:=\Tm<\infty$ and that for some $M_u>0$ \eqref{eq:u-linfty-bound} is satisfied, i.e.
\begin{align*}
\|u(\cdot,t)\|_{\Lo[\infty]}\leq M_u \quad\text{for all }t\in(0,T).
\end{align*}

In a first step we improve upon the information on $v$ by establishing a bound for the gradient in $\Lo[\infty]$.

\begin{lemma}\label{lem:nab-v-linfty}
Assume $(u_0,v_0)$ satisfy \eqref{eq:IR} with $K_{v_0}:=\|v_0\|_{\W[1,\infty]}$. Suppose that $(u,v)$ solves \eqref{eq:CT-equation} classically in $\Omega\times(0,T)$ with $T<\infty$ and that $u$ is bounded with $M_u>0$ such that \eqref{eq:u-linfty-bound} holds. Then, there is $C=C(M_u,K_{v_0})>0$ such that
\begin{align*}
\|\nabla v(\cdot,t)\|_{\Lo[\infty]}\leq C\quad\text{for all }t\in(0,T).
\end{align*}
\end{lemma}

\begin{bew}
Denoting by $\big(e^{\tau \Delta}\big)_{\tau\geq 0}$ the Neumann heat semigroup on $\Omega$, we make use of the variation-of-constants representation of $v$ to find that
\begin{align*}
\big\|\nabla v(\cdot,t)\big\|_{\Lo[\infty]}\leq \big\|\nabla e^{t\Delta}v_0\big\|_{\Lo[\infty]}+\int_0^t\big\|\nabla e^{-(t-s)\Delta} v(\cdot,t)u(\cdot,t)\big\|_{\Lo[\infty]}
\end{align*}
for all $t\in(0,T)$. The known smoothing properties of $\big(e^{\tau\Delta}\big)_{\tau\geq0}$ (see e.g. \cite[Lemma 1.3]{win10jde} and \cite[Lemma 2.1]{caolan16_smalldatasol3dnavstokes}) entail the existence of $C_1>0$ such that
\begin{align*}
\big\|\nabla v(\cdot,t)\big\|_{\Lo[\infty]}\leq C_1\big\|\nabla v_0\big\|_{\Lo[\infty]}+C_1\int_0^t\big(1+(t-s)^{-\frac{1}{2}}\big)e^{-\lambda_1 t}\big\|v(\cdot,t)u(\cdot,t)\big\|_{\Lo[\infty]}
\end{align*}
for all $t\in(0,T)$, where $\lambda_1>0$ denotes the first nonzero eigenvalue of $-\Delta$ in $\Omega$ under Neumann boundary conditions. Drawing on Lemma~\ref{lem:v-linfty_bound} and \eqref{eq:u-linfty-bound}, we hence obtain
\begin{align*}
\big\|\nabla v(\cdot,t)\big\|_{\Lo[\infty]}\leq C_1 \|\nabla v_0\|_{\Lo[\infty]}+C_1 M_u \|v_0\|_{\Lo[\infty]}\int_0^\infty \big(1+\sigma^{-\frac12}\big)e^{-\lambda_1 \sigma}\intd\sigma
\end{align*} 
for all $t\in(0,T)$. Taking $C(M_u,K_{v_0}):=C_1 K_{v_0}+C_1 M_u K_{v_0}\frac{1+\sqrt{\lambda_1\pi}}{\lambda_1}>0$ completes the proof.
\end{bew}

In the second improvement step we will turn our attention to obtaining a Hölder bound for $u$ on $\bomega\times[0,T]$ under the assumed boundedness of $u$. This will be the crucial ingredient for the Schauder theory employed in Lemma~\ref{lem:v-2-hoelder}. 

\begin{lemma}\label{lem:u-hoelder}
Assume $(u_0,v_0)$ satisfy \eqref{eq:IR} with $K_{u_0}:=\|u_0\|_{\W[1,\infty]}$ and $K_{v_0}:=\|v_0\|_{\W[1,\infty]}$. Let $M_u>0$, then there are $\theta\in(0,1)$ and $C=C(M_u,K_{u_0},K_{v_0})>0$ such that if $(u,v)$ solves \eqref{eq:CT-equation} classically in $\Omega\times(0,T)$ with $T<\infty$ and satisfies \eqref{eq:u-linfty-bound}, then
\begin{align}\label{eq:u-hoelder}
\|u\|_{\CSp{\theta,\frac{\theta}{2}}{\bomega\times[0,T]}}\leq C.
\end{align}
\end{lemma}

\begin{bew}
Writing $\Phi(s)=\int_0^s D(\sigma)\intd\sigma$, $a(x,t)=S(u)\nabla v$ and $b(x,t)=f(u,v)$ we find that if $u$ is a solution of the first equation of \eqref{eq:CT-equation} in $\Omega\times(0,T)$, it also solves
$$u_t=\Delta\Phi(u)+\nabla\cdot\big(a(x,t) u\big)+b(x,t)\quad\text{in }\Omega\times(0,T).$$
Recalling the properties of $D$ stated in \eqref{eq:cond-D}, we conclude that $\Phi\in\CSp{0}{[0,\infty)}\cap\CSp{2}{(0,\infty)}$ is convex with $\Phi(0)=0$ and $\Phi'>0$ on $(0,\infty)$ and that there are $s_0>0$ and $C_1>0$ such that
\begin{align*}
s\Phi''(s)\leq C_1 \Phi'(s)\ \text{on }[0,s_0]
\end{align*}
and that
\begin{align*}
\big\|\Phi\big(u(\cdot,t)\big)\big\|_{\Lo[\infty]}\leq \Phi(M_u)\ \text{ for all }t\in(0,T).
\end{align*}
Moreover, in view of \eqref{eq:cond-S_f} and Lemmas~\ref{lem:nab-v-linfty} and \ref{lem:v-linfty_bound}, we find $C_2=C_2(M_u,K_{v_0})>0$ such that $$\|S(u)\nabla v\|_{\LSp{\infty}{\Omega\times(0,T)}}\leq C_S(M_u)C_2\quad\text{and}\quad\|f(u,v)\|_{\LSp{\infty}{\Omega\times(0,T)}}\leq C_f(M_u,K_{v_0}).$$
Accordingly, setting $C_3=C_3(M_u,K_{v_0}):=C_S(M_u)C_2+C_f(M_u,K_{v_0})>0$ we have
\begin{align*}
\|a\|_{\LSp{\infty}{\Omega\times(0,T);\R^\dimN}}+\|b\|_{\LSp{\infty}{\Omega\times(0,T)}}\leq C_2
\end{align*}
and from $\|u_0\|_{\W[1,\infty]}\leq K_{u_0}$ and the Sobolev embedding theorem we infer the existence of some $\beta_0\in(0,1)$ such that $u_0\in\CSp{\beta_0}{\bomega}$. Since the second assumption in \eqref{eq:cond-D-hoelder} additionally ensures that $\Phi^{-1}$ is Hölder continuous on $\big[0,\Phi(M_u)\big]$ (\cite[Remark 1.4]{TB23_hoeldertaxis}), the conditions of \cite[Theorem 1.6 and Corollary 1.7]{TB23_hoeldertaxis} are satisfied and drawing on these results we obtain $\theta\in(0,1)$ and $C=C(M_u,K_{u_0},K_{v_0})>0$ such that \eqref{eq:u-hoelder}
holds.
\end{bew}

\begin{remark}
For reasonably well-behaved cross-diffusion, e.g. $D(s)=m s^{m-1}$ with $m\in(1,3]$ and $S(u)=C_S$, the lemma above can also be proved using the well-known Hölder regularity result of Porzio--Vespri (\cite{PorzVesp93}). 
\end{remark}

For the final step of our improvement procedure, we can now draw on well-known parabolic Schauder-theory to obtain a $\CSp{2+\theta,1+\frac{\theta}{2}}{\bomega\times[\tfrac{T}{2},T]}$-bound for $v$, which, of course, immediately entails the desired bound on $\|\Delta v(\cdot,t)\|_{\Lo[\infty]}$ on $(\tfrac{T}{2},T)$.

\begin{lemma}\label{lem:v-2-hoelder}
Assume that $(u_0,v_0)$ satisfy \eqref{eq:IR} with $K_{v_0}:=\|v_0\|_{\W[1,\infty]}$ and let $M_u>0$. Then, there are $\theta\in(0,1)$ and $K_1=K_1(T,M_u,K_{v_0})>0$ such that if $(u,v)$ solves \eqref{eq:CT-equation} classically in $\Omega\times(0,T)$ with $T<\infty$ and satisfies \eqref{eq:u-linfty-bound}, then
\begin{align}\label{eq:v-2-hoelder}
\|v\|_{\CSpnl{2+\theta,1+\frac{\theta}{2}}{\bomega\times[\frac{T}{2},T]}}\leq K_1.
\end{align}
In particular, there is $K_2=K_2(T,M_u,K_{v_0})>0$ such that
\begin{align}\label{eq:laplace-v-bound}
\|\Delta v(\cdot,t)\|_{\Lo[\infty]}\leq K_2\quad\text{for all }t\in\big(\tfrac{T}{2},T\big).
\end{align}
\end{lemma}

\begin{bew}
Using similar arguments as in the previous lemma, we first note that by writing $b(x,t)=-uv$ we conclude that $v$ solves $$v_t=\Delta v+b(x,t)\quad\text{in }\Omega\times(0,T)$$ with $\|b\|_{\LSp{\infty}{\Omega\times[0,T]}}\leq M_u K_{v_0}$. Since Lemma~\ref{lem:v-linfty_bound} entails that $\|v\|_{\LSp{\infty}{\Omega\times[0,T]}}\leq K_{v_0}$, and $\|v_0\|_{\W[1,\infty]}\leq K_{v_0}$ implies $v_0\in\CSp{\beta}{\bomega}$ for some $\beta\in(0,1)$, we may again draw on \cite[Theorem 1.6 and Corollary 1.7]{TB23_hoeldertaxis} or \cite[Theorem 1.3]{PorzVesp93} to find $\theta_1\in(0,1)$ and $C_1=C_1(M_u,K_{v_0})>0$ such that 
\begin{align}\label{eq:v-hoelder}
\|v\|_{\CSp{\theta_1,\frac{\theta_1}{2}}{\bomega\times[0,T]}}\leq C_1.
\end{align}
Next, we pick a smooth and monotonically increasing function $\chi:[0,T]\to\R$, satisfying $\chi\equiv 0$ on $[0,\tfrac{T}{4}]$, $\chi\equiv 1$ on $[\tfrac{T}{2},T]$ and $\|\chi\|_{\CSpnl{1}{(\frac{T}{4},\frac{T}{2})}}\leq 1+\frac{8}{T}$ and set $\widetilde{v}:=\chi v$. Then, we note that $\widetilde{v}$ satisfies $\widetilde{v}(\frac{T}{4})=0$, $\frac{\partial\widetilde{v}}{\partial\nu}\vert_{\romega}=0$ as well as
\begin{align*}
\mathscr{L}\,\widetilde{v}=\chi_t v\quad\text{on }\big(\tfrac{T}{4},T\big)
\end{align*}
with $\mathscr{L}\,\widetilde{v}:=\frac{\partial\widetilde{v}}{\partial t}-\Delta \widetilde{v}+ u\widetilde{v}$. In view of Lemma~\ref{lem:u-hoelder} there is $\theta_2\in(0,1)$ such that $\mathscr{L}$ is a linear parabolic differential operator with coefficients of class $\CSp{\theta_2,\frac{\theta_2}{2}}{\bomega\times[\tfrac{T}{4},T]}$. By the assumptions on $\chi$ and \eqref{eq:v-hoelder} there are $\theta_3\in(0,1)$ and $C_2=C_2(T,M_u,K_{v_0})>0$ such that $$\|\chi_tv\|_{\CSp{\theta_3,\frac{\theta_3}{2}}{\bomega\times[\frac{T}{4},T]}}\leq C_2$$
and from parabolic Schauder theory (e.g. \cite[Theorem III.5.1 and Theorem IV.5.3]{LSU}) we infer the existence of $C_3=C_3(T,M_u,K_{v_0})>0$ such that
$$\|\chi v\|_{\CSp{2+\theta_3,1+\frac{\theta_3}{2}}{\bomega\times[\frac{T}{4},T]}}\leq C_3.$$
Since $\chi\equiv 1$ on $[\tfrac{T}{2},T]$, this bound immediately entails the asserted bound in \eqref{eq:v-2-hoelder} from which \eqref{eq:laplace-v-bound} is an evident consequence.
\end{bew}

\setcounter{equation}{0} 
\section{Absence of dead-core formations before the blow-up time}\label{sec4}
In the final section, we are first going to employ a comparison argument with a suitably chosen spatially homogeneous lower solution $\underline{u}$ of the form $\underline{u}(x,t)=Ae^{-Bt}$ with sufficiently small $A>0$ (depending on $T$ and $\delta_0$) and suitably large $B>0$ (depending on $K_2$ from Lemma~\ref{lem:v-2-hoelder} as well as $C_S$ and $M_u$). Since we assume $T<\infty$, we can then establish $\delta_u$ such that Proposition~\ref{prop:1} holds true. After that, we can quite easily verify Theorem~\ref{theo:ex-alt} by a contradiction argument.

\begin{proof}[\textbf{Proof of Proposition~\ref{prop:1}}:]
Given any classical solution $(u,v)$ of \eqref{eq:CT-equation} in $\Omega\times(0,T)$, we conclude from the uniqueness result present in Lemma~\ref{lem:exist} and the extensibility criterion \eqref{eq:ex-alt-2} in Lemma~\ref{lem:v-linfty_bound}, that $m_u(\tfrac{T}{2}):=\min_{\bomega\times[0,\frac{T}{2}]} u>0$ is a well-defined positive number. In particular, we can pick some $A=A(T,\delta_0)>0$ such that $$A\leq \min\left\{\delta_0,m_u(\tfrac{T}{2})\right\}.$$
Denoting by $C_S=C_S(M_u)>0$ the constant from \eqref{eq:cond-S_f} and by $K_2=K_2(T,M_u,K_{v_0})>0$ the constant obtained in Lemma~\ref{lem:v-2-hoelder} satisfying
\begin{align}\label{eq:prop1-eq1}
S(\sigma)\leq C_S\ \text{for all }\ \sigma\in[0,M_u]\quad\text{and}\quad \|\Delta v(\cdot,t)\|_{\Lo[\infty]}\leq K_2\quad\text{for all }(\tfrac{T}{2},T),
\end{align}
respectively, we set $B=B(T,M_u,K_{v_0}):=K_2 C_S$ and introduce the spatially homogeneous function $\underline{u}(x,t):=A e^{-Bt}$. Then, we conclude from $f\geq 0$, \eqref{eq:prop1-eq1} and our choice for $B$ that
\begin{align*}
\underline{u}_t-\nabla\cdot\big(D(\underline{u})\nabla\underline{u}-\underline{u}S(\underline{u})\nabla v\big)-f(\underline{u},v)&=-B\underline{u}+\underline{u}S(\underline{u})\Delta v-f(\underline{u},v)\\
&\leq (C_S K_2-B)\underline{u}= 0\quad\text{in }\Omega\times(\tfrac{T}{2},T).
\end{align*}
Combined with the fact that we moreover have $\underline{u}(\cdot,\tfrac{T}{2})\leq A\leq u(\cdot,\tfrac{T}{2})$ in $\bomega$ and $\partial_\nu \underline{u}=\partial_\nu u=0$ on $\romega\times(\tfrac{T}{2},T)$, we may infer
$$A e^{-B t}=\underline{u}(x,t)\leq u(x,t)\quad\text{ in }\Omega\times(\tfrac{T}{2},T)$$
from an application of the parabolic comparison principle. Noticing that $Ae^{-B T}\leq Ae^{-B t}\leq A$ for all $t\in(0,T)$, we let $\delta_u=\delta_u(T,\delta_0,M_u,K_{v_0}):= Ae^{-B T}>0$ to obtain $$u(x,t)\geq \delta_u\quad\text{for all }(x,t)\in\Omega\times(0,T)$$ as claimed.
\end{proof}

Assuming $u$ to be bounded we have now a lower strictly positive bound for $u$ at hand. This lower bound can now be exploited to further refine the extensibility criterion \eqref{eq:ex-alt-2} and establish Theorem~\ref{theo:ex-alt}.

\begin{proof}[\textbf{Proof of Theorem \ref{theo:ex-alt}}:]
Denote by $(u,v)$ the maximally extended classical solution of \eqref{eq:CT-equation} provided by Lemma~\ref{lem:exist}. To prove that actually \eqref{eq:ref-ex-alt} holds, let us assume for contradiction that both $$\Tm<\infty\ \text{ and }\ \|u(\cdot,t)\|_{\Lo[\infty]}\leq M_u\ \text{ for all }t\in(0,\Tm)$$ are true. In this case, Proposition~\ref{prop:1} would entail the existence of $\delta_u=\delta_u(\Tm,\delta_0,M_u,K_{v_0})>0$ satisfying 
\begin{align*}
\left\|\tfrac{1}{u(\cdot,t)}\right\|_{\Lo[\infty]}\leq \frac{1}{\delta_u}\quad\text{for all }t\in(0,\Tm)
\end{align*}
and hence
\begin{align*}
\limsup_{t\nearrow\Tm}\Big(\big\|u(\cdot,t)\big\|_{\Lo[\infty]}+\left\|\tfrac{1}{u(\cdot,t)}\right\|_{\Lo[\infty]}\Big)\leq M_u+\frac{1}{\delta_u}<\infty.
\end{align*}
In view of the extensibility criterion in Lemma~\ref{lem:v-linfty_bound} this would imply $\Tm=\infty$, which clearly contradicts the assumption $\Tm<\infty$.
Accordingly, either $\Tm=\infty$, or $\limsup_{t\nearrow\Tm}\|u(\cdot,t)\|_{\Lo[\infty]}=\infty$. 
\end{proof}

\section*{Acknowledgements}
The author acknowledges support of the {\em Deutsche~Forschungsgemeinschaft} (Project No.~462888149).

\footnotesize{
\setlength{\bibsep}{3pt plus 0.5ex}

}

\end{document}